\documentclass[11pt]{amsart}

\newtheorem{theorem}{Theorem}[section]

\newtheorem{lemma}{Lemma}[section]
\newtheorem{corollary}{Corollary}[section]

\theoremstyle{definition}

\numberwithin{equation}{section}

\newcommand{\Nat}{{\mathbb N}}

\newcommand{\Real}{{\mathbb R}}

\newcommand{\B}{{\mathcal B}}
\newcommand{\Bnat}{{\mathcal B^{nat}}}
\newcommand{\hilbert}{{\mathcal H}}
\newcommand{\kspace}{{\mathcal K}}
\newcommand{\M}{{\textbf{M}}}

\def\hlimit{\overset{\hilbert}{\rightarrow}}

\begin{document}

\title[On Nyman-Beurling Criterion]{A strengthening of the Nyman-Beurling criterion for the Riemann Hypothesis}

\author{Luis B\'{a}ez-Duarte}

\date{18 February 2002}

\email{lbaez@ccs.internet.ve}


\begin{abstract}
Let $\rho(x)=x-[x]$, $\chi=\chi_{(0,1)}$. In $L_2(0,\infty)$ consider the subspace $\B$ generated by $\{\rho_a|a\geq1\}$
where $\rho_a(x):=\rho\left(\frac{1}{ax}\right)$. By the Nyman-Beurling criterion the Riemann hypothesis is equivalent
to the statement $\chi\in\overline{\B}$. For some time it has been conjectured, and proved in this paper, that the
Riemann hypothesis is equivalent to the stronger statement that  $\chi\in\overline{\Bnat}$ where $\Bnat$ is the much
smaller subspace generated by $\{\rho_a|a\in\Nat\}$.
\end{abstract}

\maketitle

\section{Introduction}
We denote the fractional part of $x$ by  $\rho(x)=x-[x]$, and let $\chi$ stand for the characteristic function of the
interval $(0,1]$. $\mu$ denotes the M\"{o}bius function. We shall be working in the Hilbert space

$$
\hilbert:=L_2(0,\infty),
$$
where the main object of interest is the subspace of \textit{Beurling functions}, defined
as the linear hull of the family $\{\rho_a|1\leq a\in\Real\}$ with

$$
\rho_a(x):=\rho\left(\frac{1}{ax}\right).
$$
The much smaller subspace $\Bnat$ of \textit{natural Beurling functions} is generated by $\{\rho_a|a\in\Nat\}$. The
Nyman-Beurling criterion (\cite{nyman}, \cite{beurling}) states, in a slightly modified form
\cite{notes3} (the original formulation is related to $L_2(0,1)$), that the Riemann hypothesis is equivalent to the
statement that 

$$
\chi\in\overline{\B},
$$
\ \\
but it has recently been conjectured by several authors\footnote{see \cite{BR1}, \cite{IUO}, \cite{AB}, \cite{notes3},
\cite{notes1}, \cite{conrey}, \cite{frankenhuysen},
\cite{landreau}, \cite{lee}, \cite{nikolski}, \cite{vasyunin}, \cite{vasyunin2}} that this condition
could be substituted by $\chi\in\overline{\Bnat}.$ We state this as a theorem to be proved below.

\begin{theorem}\label{natconj}
The Riemann hypothesis is equivalent to the statement that

$$
\chi\in\overline{\Bnat}.
$$
\end{theorem}
\ \\
To properly gauge the strength of this theorem note this: not only is $\Bnat$ a rather
thin subspace of $\B$ but, as is easily seen, it is also true that $\overline{\B}$ is much larger than
$\overline{\Bnat}$. 

By necessity all authors have been led in one way or another to the \textit{natural approximation}
\begin{equation}\label{natapprox}
F_n:=\sum_{a=1}^n\mu(a)\rho_a,
\end{equation}
which tends to $-\chi$ both a.e. and in $L_1$ norm when restricted to $(0,1)$ (see
\cite{BR1}), but which has been shown (\cite{IUO}, \cite{AB}) to diverge in $\hilbert$. In unpublished work one has
tried to prove convergence under the Riemann hypothesis of subsequences of $\{F_n\}$ such as when $n$ is
restricted to run along the solutions of $\sum_{a=1}^n\mu(n)=0$. Another attempt by J. B. Conrey and G. Myerson
\cite{conrey} relates to a mollification of $F_n$, the\textit{ Selberg approximation}, defined in \cite{notes3} by

$$
S_n:=\sum_{a=1}^n\mu(a)\left(1-\frac{\log a}{\log n}\right)\rho_a.
$$
\ \\  
A common problem to these sequences is that if they converge at all to $-\chi$ in $\hilbert$ they must do so very
slowly: it is known \cite{notes3} that for any $F=\sum_{k=1}^n c_k \rho_{a_k}$, $a_k\geq1$, if $N=\max a_k$, then

\begin{equation}\label{slow}
\|F-\chi\|_\hilbert\geq \frac{C}{\sqrt{\log N}},
\end{equation}
\ \\
for an absolute constant C that has recently been sharpened by J. F. Burnol \cite{burnol}. This, as well as
considerations of
\textit{ summability} of series, led the author in
\cite{AB} as well as here to try to employ symultaneously, as it were, the whole range of $a\in[1,\infty)$. Thus we
define for complex
$s$ and
$x>0$ the functions

\begin{equation}\label{fs}
f_s(x):=\sum_{a=1}^\infty \frac{\mu(a)}{a^s}\rho_a(x).
\end{equation}
For fixed $x>0$ this is a meromorphic functions of $s$ in the complex plane since

$$
f_s(x)=\frac{1}{x\zeta(s+1)}-\sum_{a\leq1/x}\frac{\mu(a)}{a^s}\left[\frac{1}{ax}\right],
$$ 
where the finite sum on the right is an entire function; thus $f_s$ is seen to be a sort of correction of
$1/\zeta(s)$. Assuming the Riemann hypothesis we shall prove for small positive $\epsilon$ that

$$
f_\epsilon\in\overline{\Bnat},
$$
and then,\textit{ unconditionally}, that

$$
f_\epsilon\hlimit -\chi,\ \ \ (\epsilon\downarrow0),
$$
so that $\chi\in\overline{\Bnat}$.

\section{The Proof}

\subsection{Two technical lemmae}
Here $s=\sigma+i\tau$ with $\sigma$ and $\tau$ real. The well-known theorem of Littlewood (see \cite{titchmarsh}
Theorem 14.25 (A)) to the effect that under the Riemann hypothesis $\sum_{a=1}^\infty\mu(a)a^{-s}$ converges to
$1/\zeta(s)$ for $\Re(s)>1/2$ has been provided in the more general setting of $\Re(s)>\alpha$ with a precise error
term by M. Balazard and E. Saias (\cite{notes1}, Lemme 2). We quote their lemma here for the sake of convenience.

\begin{lemma}\label{lemme2}
Let $1/2\leq\alpha<1$, $\delta>0$, and $\epsilon>0$. If $\zeta(s)$ does not vanish in the half-plane $\Re(s)>\alpha$,
then for $n\geq2$ and $\alpha+\delta\leq\Re(s)\leq1$ we have

\begin{equation}
\sum_{a=1}^n\frac{\mu(a)}{a^s}=
\frac{1}{\zeta(s)}+O_{\alpha,\delta,\epsilon}\left(n^{-\delta/3}(1+|\tau|)^\epsilon\right)
\end{equation} 
\end{lemma}
\ \\
It is important to note that the next lemma is independent of the Riemann or even the Lindel\"{o}f hypothesis.

\begin{lemma}\label{zratio}
For $0\leq\epsilon\leq \epsilon_0<1/4$ there is a positive constant $C=C(\epsilon_0)$ such that
for all
$\tau$
\begin{equation}\label{zratioeq}
\left|\frac{\zeta(\frac{1}{2}-\epsilon+i\tau)}{\zeta(\frac{1}{2}+\epsilon+i\tau)}\right|\leq
C\left(1+|\tau|\right)^\epsilon.
\end{equation}
\end{lemma}
\begin{proof}
We bring in the functional equation of $\zeta(s)$ to bear as follows
\begin{eqnarray}\nonumber
\left|\frac{\zeta(\frac{1}{2}-\epsilon+i\tau)}{\zeta(\frac{1}{2}+\epsilon+i\tau)}\right|&=&
\left|\frac{\zeta(\frac{1}{2}-\epsilon-i\tau)}{\zeta(\frac{1}{2}+\epsilon+i\tau)}\right|\\\nonumber
&=&\pi^{-\epsilon}\left|
\frac{\Gamma(\frac{1}{4}+\frac{1}{2}\epsilon+\frac{1}{2}i\tau)}{\Gamma(\frac{1}{4}+\frac{1}{2}\epsilon+\frac{1}{2}i\tau)}
\right|,
\end{eqnarray}
then the conclusion easily follows from well-known asymptotic formulae for the gamma function in a vertical strip
(\cite{rademacher} (21.51), (21.52)).
\end{proof}

\subsection{The proof proper of Theorem \ref{natconj}} 

It is clear that we need not prove the \textit{if} part of Theorem \ref{natconj}. So let us assume that the Riemann
hypothesis is true. We define

$$
f_{\epsilon,n}:=\sum_{a=1}^n\frac{\mu(a)}{a^\epsilon}\rho_a,\ \ \ (\epsilon>0).
$$ 
It is easy to see that

\begin{equation}\label{fepsilonn}
f_{\epsilon,n}(x)=\frac{1}{x}\sum_{a=1}^n\frac{\mu(a)}{a^{1+\epsilon}}
-\sum_{a=1}^n\frac{\mu(a)}{a^\epsilon}\left[\frac{1}{ax}\right],
\end{equation}
then, noting that the terms of the right-hand sum drop out when $a>1/x$, we obtain the pointwise limit

\begin{equation}\label{fepsilon1}
f_\epsilon(x)=\lim_{n\rightarrow\infty}f_{\epsilon,n}(x)=\frac{1}{x\zeta(1+\epsilon)}
-\sum_{a\leq 1/x}\frac{\mu(a)}{a^\epsilon}\left[\frac{1}{ax}\right].
\end{equation}
Then again for fixed $x>0$ we have
\begin{equation}\label{endlimit}
\lim_{\epsilon\downarrow0}f_\epsilon(x)=-\sum_{a\leq 1/x}\mu(a)\left[\frac{1}{ax}\right]=-\chi(x),
\end{equation}
by the fundamental property on M\"{o}bius numbers. The task at hand now is to prove these pointwise limits are also
valid in the $\hilbert$-norm. To this effect we introduce a new Hilbert space

$$
\kspace:=L_2((\infty,\infty),(2\pi)^{-1/2}dt),
$$
\ \\
and note that by virtue of Plancherel's theorem the Fourier-Mellin map
$\M$ defined by

\begin{equation}\label{isometry}
\M(f)(\tau):=\int_0^\infty x^{-\frac{1}{2}+i\tau}f(x)dx,
\end{equation}
\ \\
is an invertible isometry from $\hilbert$ to $\kspace$. A well-known identity, which is at the root of the
Nyman-Beurling formulation, probably due to Titchmarsh (\cite{titchmarsh}, (2.1.5)), namely

$$
-\frac{\zeta(s)}{s}=\int_0^\infty x^{s-1}\rho_1(x)dx,\ \ \ (0<\Re(s)<1),
$$
immediately yields, denoting $X_{\epsilon}(x)=x^{-\epsilon}$,

\begin{equation}\label{mfepsilonn}
\M(X_\epsilon f_{2\epsilon,n})(\tau)=
-\frac{\zeta(\frac{1}{2}-\epsilon+i\tau)}{\frac{1}{2}-\epsilon+i\tau}
\sum_{a=1}^n \frac{\mu(a)}{a^{\frac{1}{2}+\epsilon+i\tau}},\ \ \ (0<\epsilon<1/2).
\end{equation}
\ \\
By a theorem of Littlewood (\cite{titchmarsh}, Theorem 14.25 (A)) if we let $n\rightarrow\infty$ in the right-hand side
of (\ref{mfepsilonn}) we get the pointwise limit

\begin{equation}\label{prelimit}
-\frac{\zeta(\frac{1}{2}-\epsilon+i\tau)}{\frac{1}{2}-\epsilon+i\tau}
\sum_{a=1}^n \frac{\mu(a)}{a^{\frac{1}{2}+\epsilon+i\tau}} \rightarrow
-\frac{\zeta(\frac{1}{2}-\epsilon+i\tau)}{\zeta(\frac{1}{2}+\epsilon+i\tau)}\frac{1}{\frac{1}{2}-\epsilon+i\tau}.
\end{equation}
\ \\
To see that this limit also takes place in $\hilbert$ we choose the parameters in Lemma \ref{lemme2} as
$\alpha=1/2$, $\delta=\epsilon>0$, $\epsilon\leq 1/2$, and $n\geq2$ to obtain

$$
\sum_{a=1}^n \frac{\mu(a)}{a^{\frac{1}{2}+\epsilon+i\tau}}=
\frac{1}{\zeta(\frac{1}{2}+\epsilon+i\tau)}+O_\epsilon\left((1+|\tau|)^\epsilon\right).
$$ 
If we now use Lemma \ref{zratio} and the Lindel\"{o}f hypothesis applied to the abcissa $1/2-\epsilon$, which follows
from the Riemann hypothesis, we obtain a positive constant $K_\epsilon$ such that for all real $\tau$

$$
\left|-\frac{\zeta(\frac{1}{2}-\epsilon+i\tau)}{\frac{1}{2}-\epsilon+i\tau}
\sum_{a=1}^n \frac{\mu(a)}{a^{\frac{1}{2}+\epsilon+i\tau}}\right|
\leq K_\epsilon(1+|\tau|)^{-1+2\epsilon}.
$$
\ \\
It is then clear that for $0<\epsilon<1/4$ the left-hand side of (\ref{prelimit}) is uniformly majorized by a
function in $\kspace$. Thus the convergence does take place in $\kspace$ which implies that

$$
X_\epsilon f_{2\epsilon,n}\hlimit X_\epsilon f_{2\epsilon}.
$$
\ \\
But $x^{-\epsilon}>1$ for $0<x<1$, and for $x>1$

\begin{equation}\label{forx>1}
f_{2\epsilon,n}(x)=\frac{1}{x}\sum_{a=1}^n\frac{\mu(a)}{a^{1+\epsilon}}\ll\frac{1}{x}, \ \ \ (x>1),
\end{equation}
which easily implies that one also has $\hilbert$-convergence for $f_{2\epsilon,n}$ as $n\rightarrow\infty$. The
factor $2$ in the subindex is unessential, so that we now have for sufficiently small $\epsilon>0$ that

$$
f_{\epsilon,n}\hlimit f_\epsilon\in\overline{\Bnat},
$$
\ \\  
as was announced above. Moreover, since we have identified the pointwise limit in (\ref{prelimit}) we now
have

$$
\M(X_\epsilon f_{2\epsilon})(t)=
-\frac{\zeta(\frac{1}{2}-\epsilon+i\tau)}{\zeta(\frac{1}{2}+\epsilon+i\tau)}\frac{1}{\frac{1}{2}-\epsilon+i\tau}.
$$
\ \\
Now we apply Lemma \ref{zratio} and obtain, \textit{without the assumption of the Riemann hypothesis}, that
$\M(X_\epsilon f_{2\epsilon})$ converges in $\kspace$, thus $X_\epsilon f_{2\epsilon}$ converges in $\hilbert$, and
this means that $f_\epsilon$ also converges in $\hilbert$ as $\epsilon\downarrow0$ by an argument entirely similar to
that used for $f_{\epsilon,n}$. The identification of the pointwise limit in (\ref{endlimit}) finally gives

$$
f_\epsilon\hlimit -\chi,
$$   
\ \\
which concludes the proof.

\section{Some comments and a corollary}
The proof of Theorem \ref{natconj} provides in turn a new proof, albeit of a stronger theorem, of the
Nyman-Beurling criterion which bypasses the deep and complicated Hardy space techniques. One should extend it to the
$L_p$ case, that is, to the condition that $\zeta(s)$ does not vanish in a half-plane $\Re(s)>1/p$. It should be clear
also that we have shown this special equivalence criterion to be true:
\begin{corollary}
The Riemann hypothesis is equivalent to the $\hilbert$-convergence of $f_{\epsilon,n}$ as
$n\rightarrow\infty$ for all sufficiently small $\epsilon>0$. 
\end{corollary}
In essence what has been done is to apply a summability method to the old natural approximation. The convergence on
special subsequences both of $n$ and of $\epsilon$ is also necessary and sufficient, and it is proposed here to study
this alongside with other summability methods for the natural approximation.
\ \\
A final remark is in order. Note that we did not employ the dependence on $n$ in the Balazard-Saias Lemma
\ref{lemme2}. This dependence would seem to be closely connected to the slowness of approximation to $-\chi$ indicated
in (\ref{slow}). 
\ \\
\ \\
\textbf{Ackowledgements.} The author wants to thank M. Balazard and E. Saias for pointing out their important, delicate
lemmae \ref{lemme2}.

\bibliographystyle{amsplain}

\ \\
\noindent Luis B\'{a}ez-Duarte\\
Departamento de Matem\'{a}ticas\\
Instituto Venezolano de Investigaciones Cient\'{\i}ficas\\
Apartado 21827, Caracas 1020-A\\
Venezuela\\

\end{document}